\documentclass[11pt]{article}

\usepackage{bridges}
\usepackage{amsfonts,amssymb,amsthm,eucal,amsmath}
\usepackage{graphicx}
\usepackage{subfig}

\usepackage{wrapfig}
\usepackage{pinlabel, color}
\usepackage{multirow}

\usepackage{microtype}
\usepackage{hyperref}

\usepackage[para]{manyfoot}
\DeclareNewFootnote[para]{B}

\newcommand{\EE}{\mathbb{E}}
\newcommand{\HH}{\mathbb{H}}

\newcommand{\RR}{\mathbb{R}}

\newcommand{\SO}{\text{SO}}

\newcommand{\dr}{d\mathbf{r}}

\theoremstyle{definition}

\newtheorem*{defn*}{Definition}

\begin{document}

\title{Non-euclidean virtual reality I: explorations of $\HH^3$}
\author{
 \begin{tabular}{cccc}
  Vi Hart & Andrea Hawksley & Elisabetta A. Matsumoto & Henry Segerman \\
  eleVR & eleVR & School of Physics & Department of Mathematics  \\
  HARC & HARC & Georgia Institute of Technology & Oklahoma State University 
 \end{tabular}
}
\date{}

\maketitle

\begin{abstract}
We describe our initial explorations in simulating non-euclidean geometries in virtual reality. Our simulations of three-dimensional hyperbolic space are available at \href{http://h3.hypernom.com}{h3.hypernom.com}.\footnote{The code is available at \href{https://github.com/hawksley/hypVR}{github.com/hawksley/hypVR}.} 
\end{abstract}

\begin{figure}[htb]
\centering
{
\includegraphics[width=\textwidth]{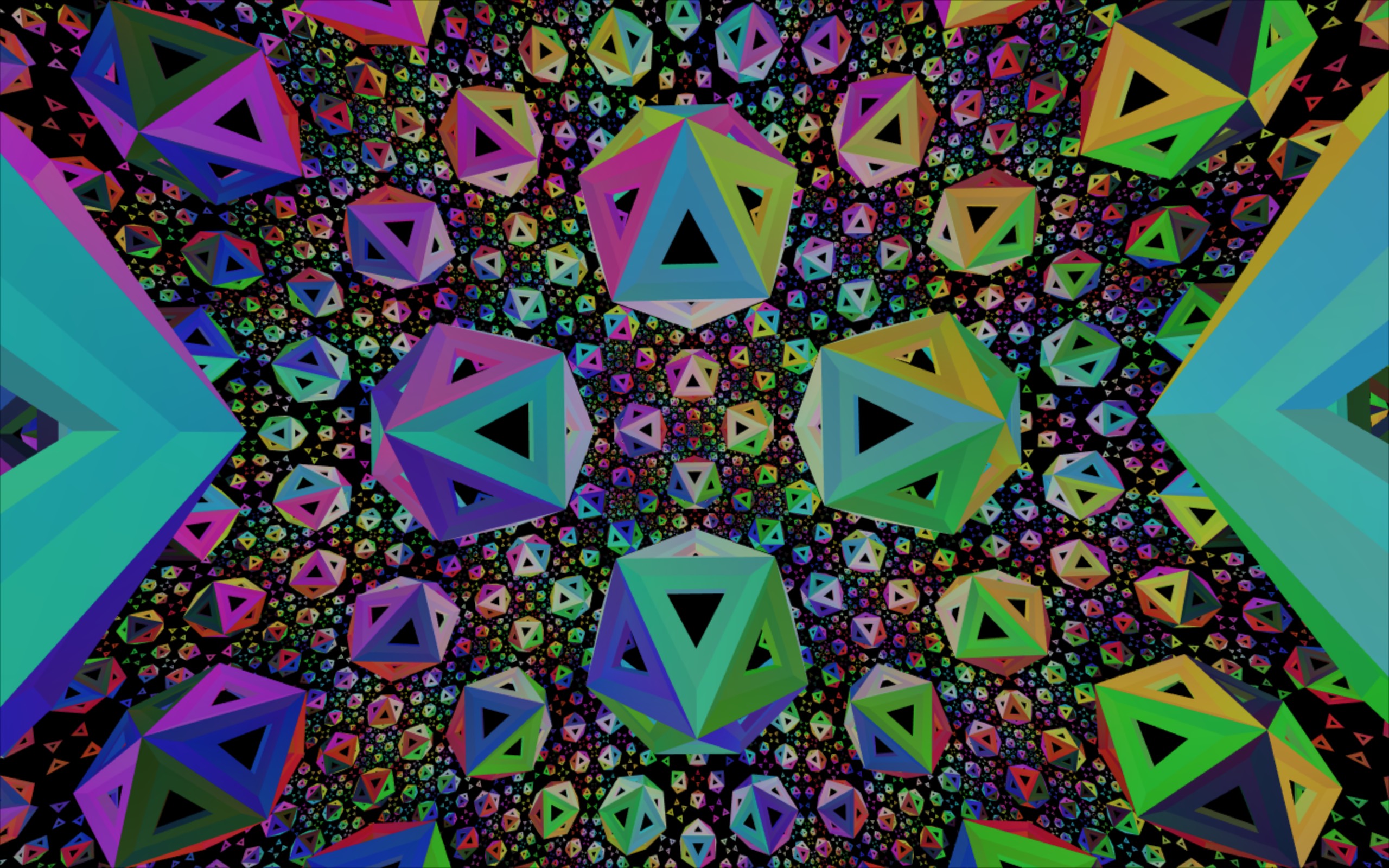}
\label{h2xr_title}
}
\caption{A view from $\HH^3$.}
\label{H2xE_title}
\end{figure}

We are used to living in three-dimensional euclidean space, and our day-to-day experiences of curvature centre around surfaces embedded in $\EE^3$. In the study of topology, the closed two-dimensional surfaces are the sphere, the torus, the two-holed torus, the three-holed torus, and so on. Thinking of these surfaces topologically, they don't come with a particular choice of \emph{geometry} -- that is, we can think about a surface as if they were made from plasticine -- without knowing what lengths and angles mean on the surface. There are however particularly nice geometries for these surfaces:  \emph{isotropic} geometries, meaning that the geometry is the same everywhere in the space, facing in every direction. We have spherical geometry for the sphere, euclidean geometry for the torus, and hyperbolic geometry for all of the others. 
In three dimensions, the story is more complicated. Thurston's geometrization conjecture, proved by Perelman~\cite{perelman1}, gives eight geometries that a three-manifold can take (although the manifold may need to be decomposed into pieces, each with one of the eight geometries). The eight geometries are $S^3$, $\EE^3$, $\HH^3$, $S^2\times\EE$, $\HH^2\times\EE$, Nil, Solv, and $\widetilde{\text{PSL}_2(\RR)}$~\cite{thurston_book}. The first three are again isotropic: spherical, three-dimensional euclidean and hyperbolic geometries. The second two are mixtures of the two-dimensional geometries and one-dimensional euclidean space, and so are not isotropic: the geometry looks different when you look in the euclidean versus the non-euclidean directions. The last three are more complicated ``twisted'' versions of these mixed geometries.

\begin{figure}[tb]
\centering
\vspace{-10pt}\subfloat[The three torus, giving the $\{4,3,4\}$ honeycomb.]
{
\includegraphics[width=0.47\textwidth]{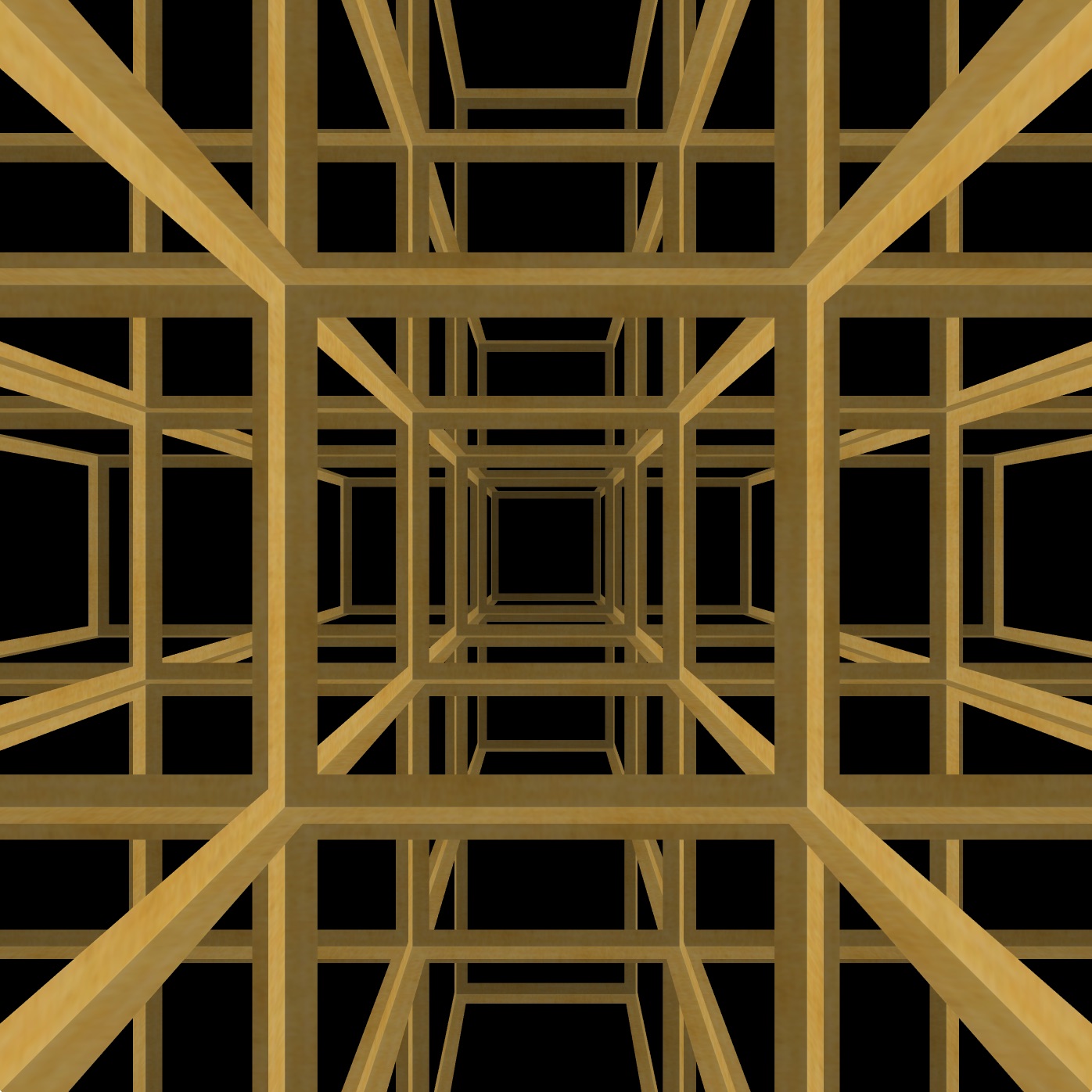}
\label{curved_spaces_cubes}
}
\quad
\subfloat[The $\{5,3,4\}$ honeycomb.]
{
\includegraphics[width=0.47\textwidth]{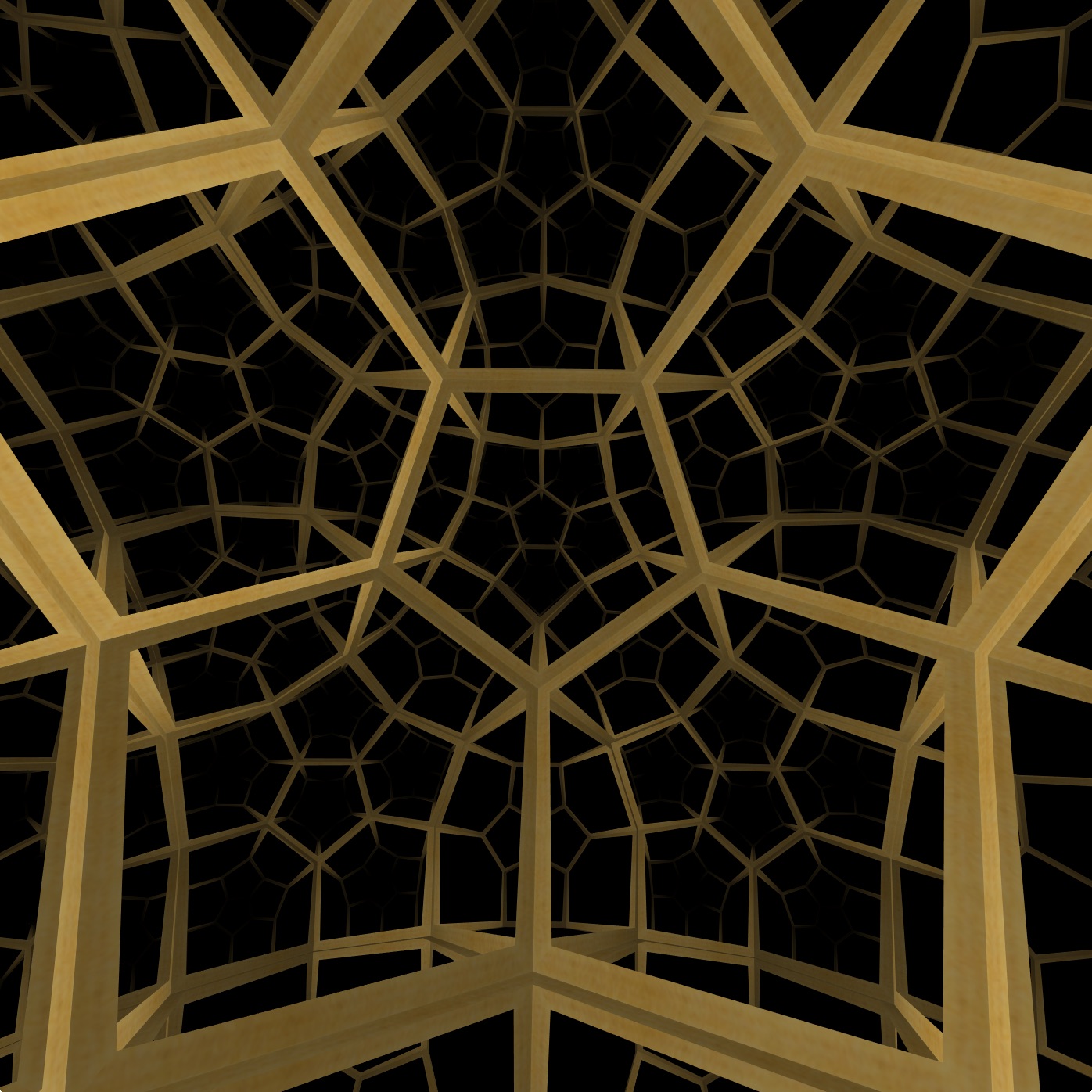}
\label{curved_spaces_dodecahedra}
}
\caption{Screenshots from \emph{Curved spaces} by Jeff Weeks.}
\label{curved_spaces}
\end{figure}

Jeff Weeks' software \emph{Curved Spaces}~\cite{curved_spaces} is a ``flight simulator for multiconnected universes''. See Figure \ref{curved_spaces}. It simulates what it would be like to explore a selection of closed three-dimensional manifolds, with $S^3$, $\EE^3$ and $\HH^3$ geometries. 
Each of these are viewed as if we are living inside the space and seeing objects in that space via rays of light that travel along \emph{geodesics} in the space.
That is, light travels along paths of shortest distance.

We are currently developing a virtual reality simulation of $\HH^3$, using many of the same ideas as are used in Weeks' work. Weeks explains the implementation in detail in \cite{weeks_real-time_rendering}; we give an overview in 
this paper.
Positional tracking in modern virtual reality headsets lets us experience
features of hyperbolic space, such as the effects of parallel transport, in a very direct way.

There are four ingredients that go into our virtual reality simulation of $\HH^3$ as outlined in this paper:
\begin{enumerate}
\item{A way to describe the points of $\HH^3$ numerically, i.e. a \emph{model} of $\HH^3$}

\item{A way to convert points in the model into points in $\EE^3$ that we can then draw on screen,}

\item{A way to move around $\HH^3$ using the motion inputs from the virtual reality headset, and}

\item{A set of landmarks in $\HH^3$ to draw, to help the viewer navigate the space -- 
we use a tiling of $\HH^3$.}

\end{enumerate}

\newpage
\section{The Model of $\HH^3$}

\begin{wrapfigure}[12]{r}{0.45\textwidth}
  \vspace{-5pt}
  \centering
  \includegraphics[width=0.4\textwidth]{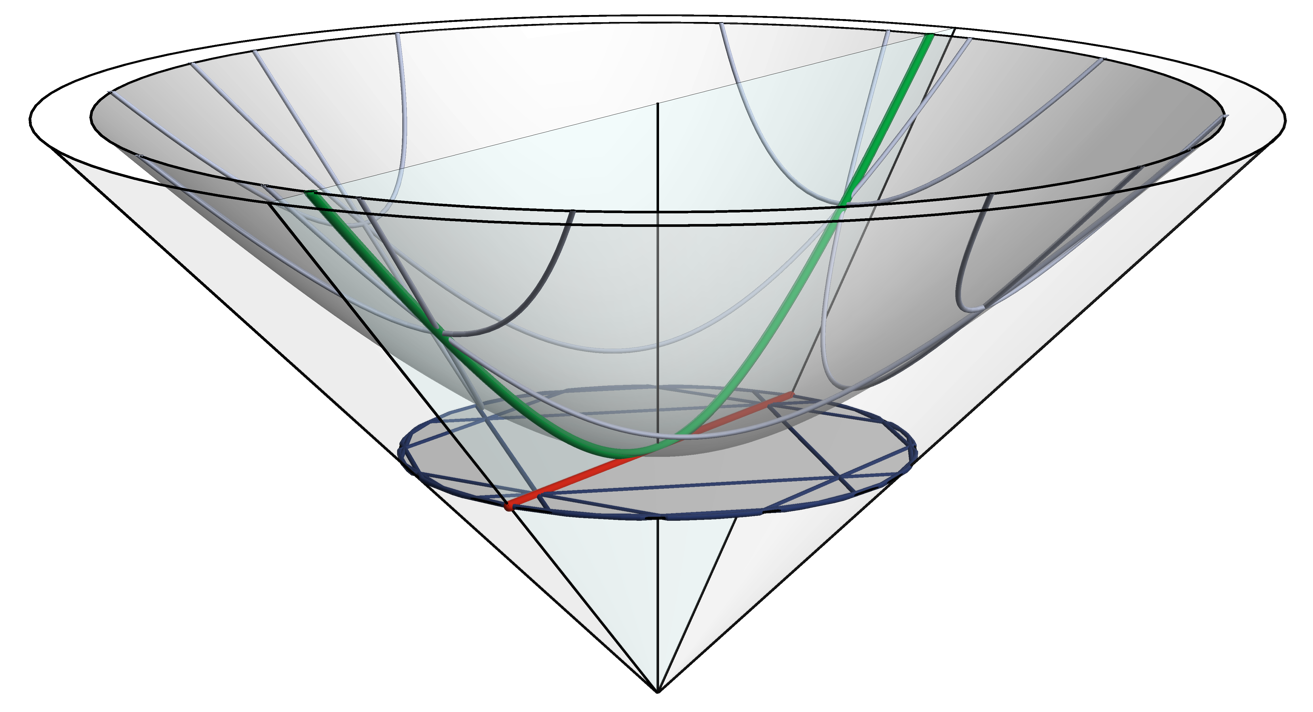}
  \caption{The hyperboloid and Klein models of $\HH^2$. Projecting the hyperboloid towards the origin onto the plane $w=1$ results in the Klein model.} 
  \label{Fig:hyperboloid}
 \end{wrapfigure}
For the first ingredient, there are many different models of hyperbolic space, including the Poincar\'e disk model, the upper half plane model, the Klein model, and the hyperboloid model. 
Compared to the other commonly seen models, the hyperboloid model is less easy to use for direct visualisation, but it turns out to be very useful for calculation. The hyperboloid model of $\HH^2$ is the set of points ${\{(x,y,w) \in \EE^{2,1} \mid x^2+y^2 = w^2 - 1, w>0\}}$, where $\EE^{2,1}$ is Minkowski space with two space-like directions, $x,y$ and one time-like direction, $w$.   Three-dimensional Minkowski space $\EE^{2,1}$ has the same cartesian coordinate system as $\EE^3$, but comes equipped with a different metric, which has line element $ds^2=dx^2+dy^2-dw^2$. The \emph{metric} $g_{ij}$ is a tensorial function that generalises the method of computing distances and angles (i.e. the dot product in euclidean space) to a differentiable manifold.

With the metric induced from the Minkowski space it lives in, the hyperboloid then has constant gaussian curvature $-1$, i.e. it is a model for the hyperbolic plane. 
For each point $(x,y,w)$ of the hyperboloid, we can divide the coordinates by $w$ to obtain $(x/w,y/w,1)$. This maps the hyperboloid to the unit radius disk on the $w=1$ plane.
The result is the Klein model of $\HH^2$.
See Figure \ref{Fig:hyperboloid}. 
Geodesics in the hyperboloid model of $\HH^2$ are intersections of the hyperboloid with planes in $\EE^{2,1}$ that pass through $(0,0,0)$. 
These geodesics map to the Klein model of $\HH^2$ as straight lines (in the euclidean sense).

Three-dimensional hyperbolic space $\HH^3$, and indeed higher dimensional hyperbolic spaces, can be modeled analogously to $\HH^2$, but in higher dimensional ambient spaces. The generalised hyperboloid model for $d-$dimensional hyperbolic space $\HH^d$ 
is the set of points in Minkowski space with $d$ space-like directions, $x_1,...,x_d$ and one time-like direction $w$, given by $\{(x_1,x_2,..,x_d,w) \in \EE^{d,1} \mid \sum_{n=1}^{d} x_n^2 = w^2 - 1, w>0\}.$

\section{Drawing points in $\HH^3$ on screen}

In order to draw a point of $\HH^3$ on the screen, we need to understand the relationship between the location of the point on the hyperboloid and us, the viewer, situated at the \emph{origin} of the hyperboloid, $(0,0,0,1)\in\EE^{3,1}$. 
We are not actually viewing points in $\HH^3$, but we view their image in the \emph{tangent space} at the origin  -- a copy of $\EE^3$ consisting of the tangent vectors of the hyperboloid at the origin. A point $p_{\HH^3} \in \HH^3$ is connected to the origin by a parametrised geodesic $\boldsymbol{\gamma}(t)$ that leaves the origin at $t=0$ and intercepts $p_{\HH^3}$ at $t=1$.
Our view of the same point $p_{\EE^3} \in \EE^3$ should also be connected to us via a geodesic in $\EE^3$ (i.e. a straight line). The velocity of the geodesic in $\HH^3$ at the origin $\dot{\boldsymbol{\gamma}}(0),$ tells us where to find $p_{\EE^3}$ -- its direction is the direction in which we must look to find $p_{\EE^3}$, and its magnitude indicates the distance between us (situated at the origin of the hyperboloid) and $p_{\EE^3}$. The map we have described, taking points on the hyperboloid to points in $\RR^3$, is the inverse of the \emph{(riemannian geometry) exponential map}. 
The exponential map goes in the other direction, sending points in the tangent space of the hyperboloid at our location into the hyperboloid.

\begin{figure}[htb]
\centering
\subfloat[Cubes.]
{
\includegraphics[width=0.4\textwidth]{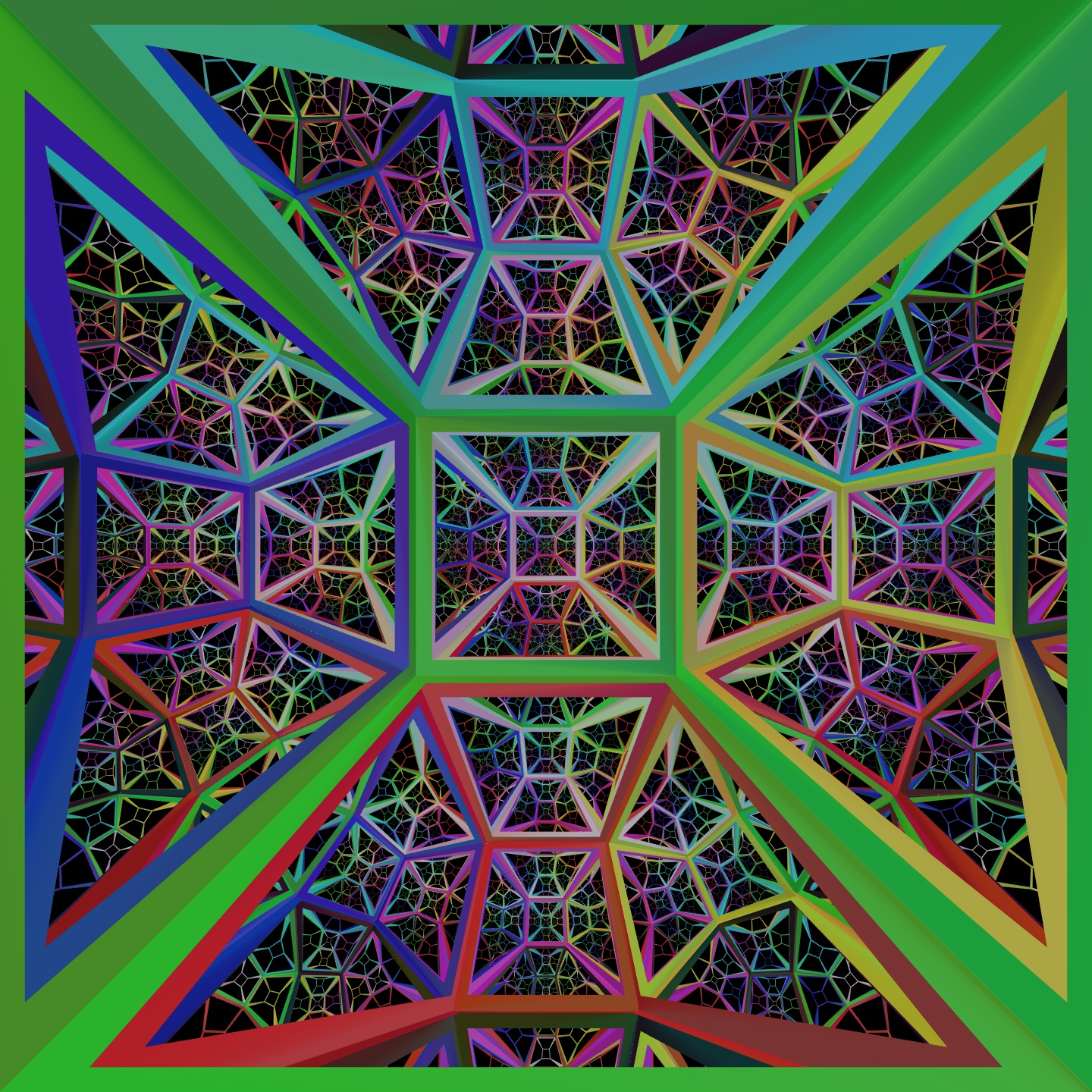}
\label{436_cubes}
}
\quad
\subfloat[Truncated cubes.]
{
\includegraphics[width=0.4\textwidth]{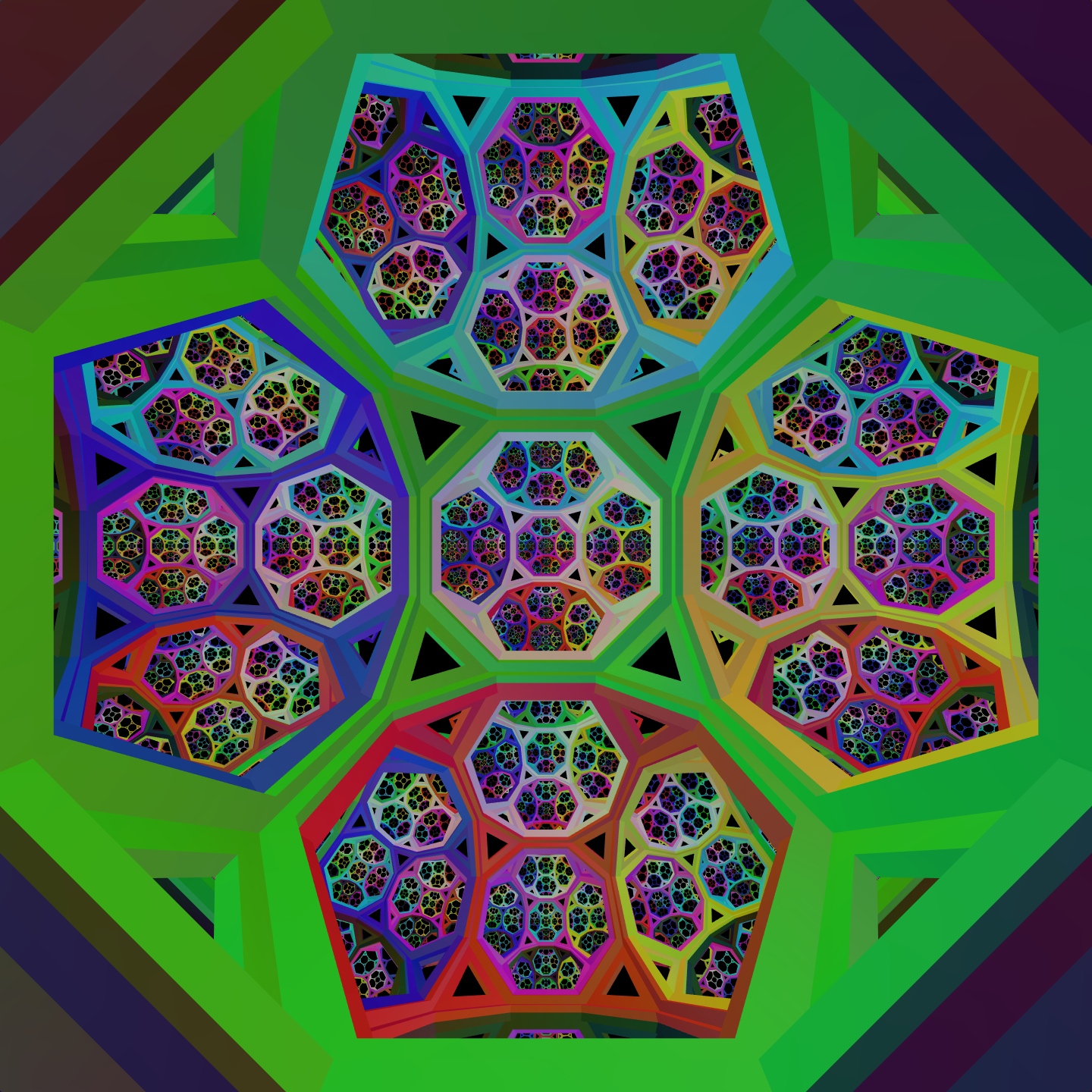}
\label{436_trunc_cubes}
}

\vspace{-5pt}
\subfloat[Only the triangular faces of the truncated cubes.]
{
\includegraphics[width=0.4\textwidth]{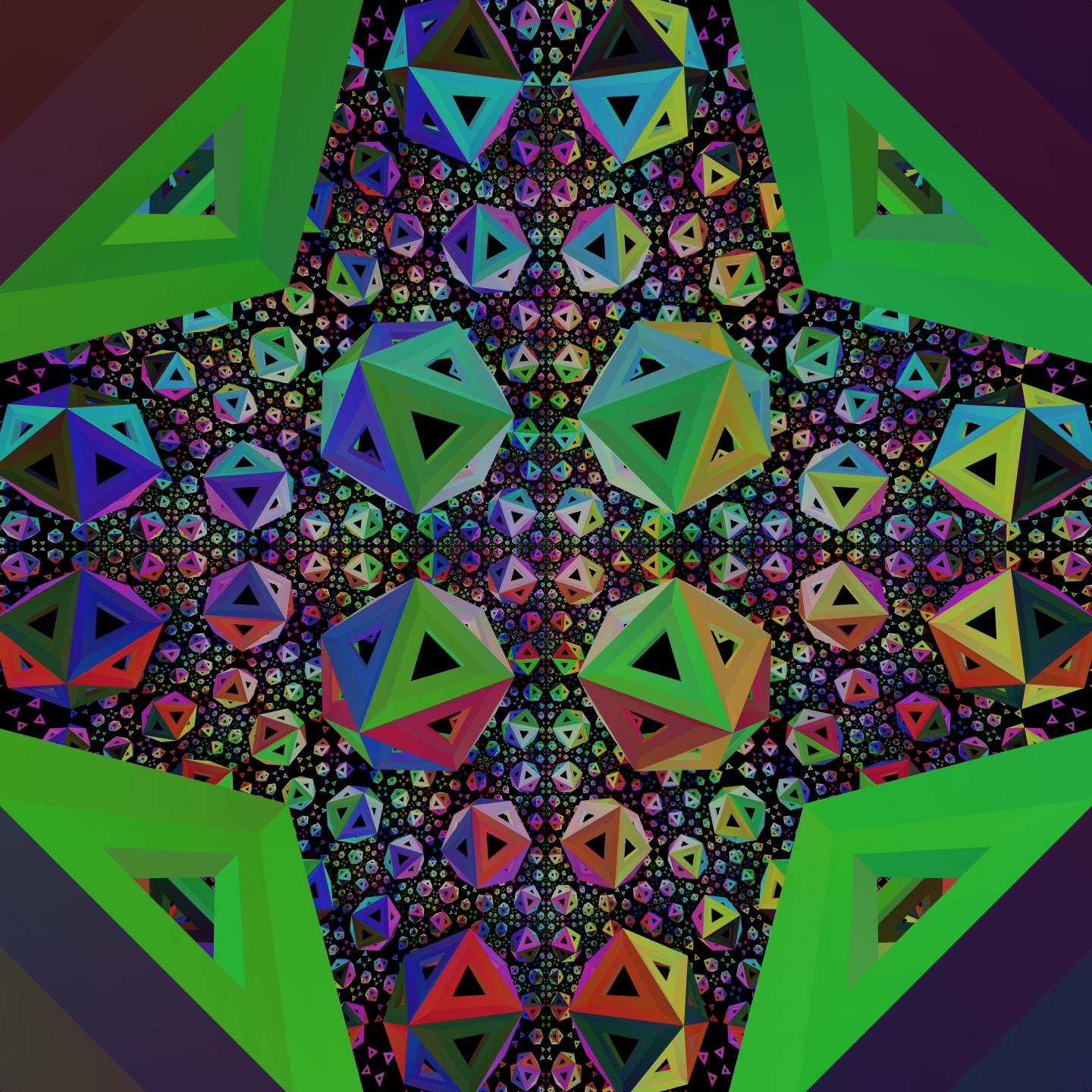}
\label{436_trunc_cubes_tris_only}
}
\quad
\subfloat[The view from inside the ``polyhedron'' in the centre lower left of Figure \ref{436_trunc_cubes_tris_only.jpg}.]
{
\includegraphics[width=0.4\textwidth]{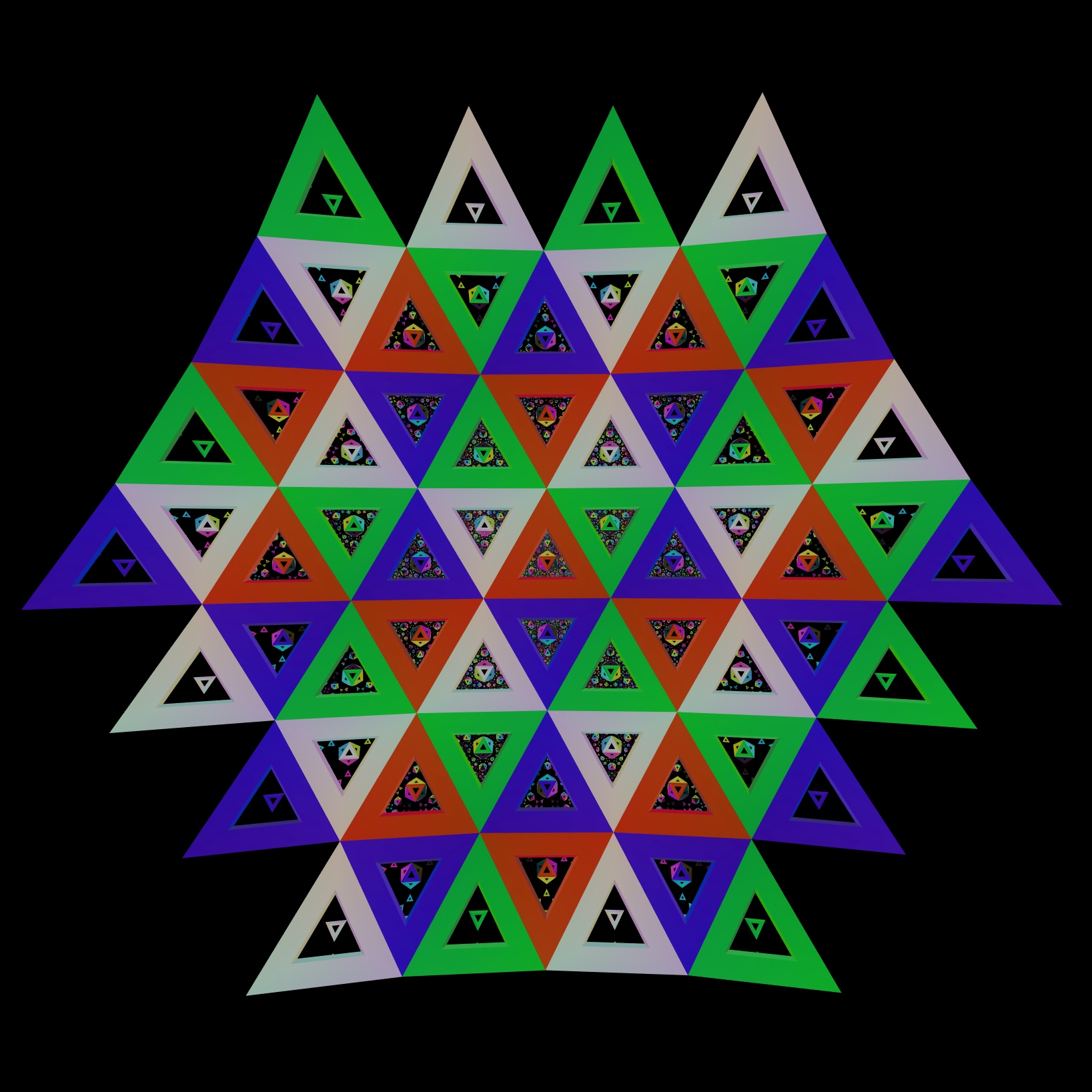}
\label{view_from_a_cusp}
}
\caption{Views of the $\{4,3,6\}$ honeycomb. We draw the honeycomb out to a depth of six steps from the central cube.}
\label{436}
\end{figure}

The correct thing to do is to use the inverse of the exponential map to draw points on screen, but because $\HH^3$ is isotropic, as are most of the models used to draw it, we can get by without calculating this. Its isotropy implies that the viewer cannot tell the difference between looking off in two different directions without aid of the decorations we use as landmarks. Likewise, we don't actually need to compute the absolute distance a point is from us. We merely need to know the relative distance between two points in a given direction, so that nearer points appear closer to us. The Klein model is computationally the cheapest to calculate -- as it does not involve inverse hyperbolic trigonometric functions -- so this is the one we (and Weeks~\cite{curved_spaces}) choose. Figure \ref{436} shows a number of views of a honeycomb in $\HH^3$ drawn using this algorithm.

\vspace{-5pt}
\section{Moving through $\HH^3$}

Although the trick of implementing graphics using the Klein model only works at the origin, we can still leverage its power as we move through the space.
As in many computer graphics implementations, we leave the viewer at the origin and translate the world around them to simulate the viewer's movement. The appropriate ``translations'' for us are isometries of $\HH^3.$ Infinitesimal translations are given by the generators of the Lie group of the space and finite transformations are given by the \emph{(Lie theory) exponential map}\footnote{This is similar to the riemannian geometry version of the exponential map, except that instead of converting a tangent vector (an infinitesimal movement in some direction) into a point at the end of a geodesic segment, it converts a more general infinitesimal motion into an isometry.}. As with $\HH^2$, the isometries of $\HH^3$ are isometries of $\EE^{3,1}$ which preserve the hyperboloid and its metric. These are elements of the group $\SO(3,1).$ The translation by a vector $\dr=(dx,dy,dz)$ in the tangent space is given by the exponential $\exp(\mathbf{M})=\sum_{n=0}^\infty \frac{1}{n!}\mathbf{M}^n$ of the matrix 
\[\mathbf{M}= \left(
\begin{array}{cccc}
0 & 0 & 0 & dx\\
0 & 0 & 0 & dy\\
0 & 0 & 0 & dz\\
dx & dy & dz & 0
\end{array}
\right).\]
Calculation of the series for the matrix exponential $\exp(\mathbf{M})$ can be vastly simplified due to a trick pointed out by Jeff Weeks. Note that $\mathbf{M}^3=|\dr|^2\mathbf{M},$ and $\mathbf{M}^4=|\dr|^2\mathbf{M}^2$, where $|\dr| = \sqrt{dx^2+dy^2+dz^2}$. Then the matrix exponential can be split into two sums: 
\[
\sum_{n=1}^\infty \frac{|\dr|^{2n-2}}{(2n-1)!}\mathbf{M}=\frac{\sinh(|\dr|)}{|\dr|}\mathbf{M},
\quad
\sum_{n=1}^\infty\frac{|\dr|^{2n}}{(2n)!}|\dr|^{2n}\mathbf{M}=\frac{\cosh(|\dr|)-1}{|\dr|^2}\mathbf{M}.
\] 
Thus, the exponential map is given by 
\[\exp\mathbf{M}=\mathbf{Id}+\frac{\sinh(|\dr|)}{|\dr|}\mathbf{M}+\frac{\cosh(|\dr|)-1}{|\dr|^2}\mathbf{M}^2.\]

When the user moves their head, the virtual reality headset detects this movement in the three-dimensional euclidean space in which we live. The difference in position between two subsequent frames is some vector, which gives us the translation of the user $-\dr$.\footnote{Note that the sign of $-\dr$ is due to the fact that the sensors detect the displacement of the virtual reality headset as $-\dr$, which corresponds to moving the entire world by a vector with the same magnitude, but in the opposite direction, $\dr$.} We then generate the isometry $\exp(\mathbf{M})$, and apply it to all the points of our simulated world before rendering the next frame. This moves the points of the world in the hyperboloid by isometries, giving the correct sense in which the user moves through the world.

\section{Decoration: the $\{4,3,6\}$ honeycomb and its colouring}

Any three-dimensional manifold can be made by taking a polyhedron and gluing its sides together in some way. Jeff Weeks' Curved Spaces shows such a polyhedron for each manifold. For example, Figure \ref{curved_spaces_cubes} shows the view from inside the three-torus, whose geometry is $\EE^3$. In this case, the polyhedron is a cube with opposite sides glued. We see a tiling (or honeycomb) of $\EE^3$ by cubes -- what we get by ``unwrapping'' the three-torus into space. This tiling has \emph{Schl\"afli symbol} $\{4,3,4\}$, meaning that the faces are squares (with 4 sides), the cells are made out of these faces, with 3 around each vertex, and there are 4 cells arranged around each edge. Figure \ref{curved_spaces_dodecahedra} shows a different manifold with a corresponding honeycomb in which four dodecahedra meet around each edge. The corresponding Schl\"afli symbol is $\{5,3,4\}$, corresponding to cells made out of pentagons (5 sides), with 3 around each vertex, and with 4 cells arranged around each edge. 

As our fourth ingredient, we decorate $\HH^3$ with another honeycomb of cubes, this time with Schl\"afli symbol $\{4,3,6\}$. See Figure \ref{436_cubes}. Here we have six cubes around each edge, rather than four. A surprising feature of this honeycomb is that the cubes are no longer of finite size -- it turns out that the vertices must be infinitely far away. See \cite{visualizing_hyperbolic_honeycombs} for more on this phenomenon.

For the euclidean honeycomb $\{4,3,4\}$, with four cubes around each edge, if we truncate each of the cubes, cutting off the corners, the revealed triangular faces form an octahedron arranged around each vertex of the original honeycomb. If we do the same thing for our hyperbolic honeycomb, as in Figure \ref{436_trunc_cubes}, the triangular faces form an infinite tiling -- the tiling of the euclidean plane with six triangles around each vertex. In our visualisation, we can experience this directly. This is easiest to see if we remove the edges of the cubes, leaving only the triangular faces, as in Figure \ref{436_trunc_cubes_tris_only}. These form strange looking polyhedra at first sight: one could believe that they are icosahedra, except that the vertex degree is six.
If you put your head ``into'' one of these polyhedra, and look back out from the inside, the polyhedron becomes the tiling of the euclidean plane, as we see in Figure \ref{view_from_a_cusp}. These polyhedra in fact correspond to \emph{horospheres} in $\HH^3$. These are ``spheres'' centered on points on the boundary of $\HH^3$, whose induced metric is the same as the euclidean plane -- which allows us to draw the regular tiling by equilateral triangles on them seen in Figure \ref{view_from_a_cusp}.

\begin{wrapfigure}[15]{l}{0.65\textwidth}
\vspace{-15pt}
\centering
\subfloat[Colouring of the hypercube, view 1.]
{
\includegraphics[width=0.3\textwidth]{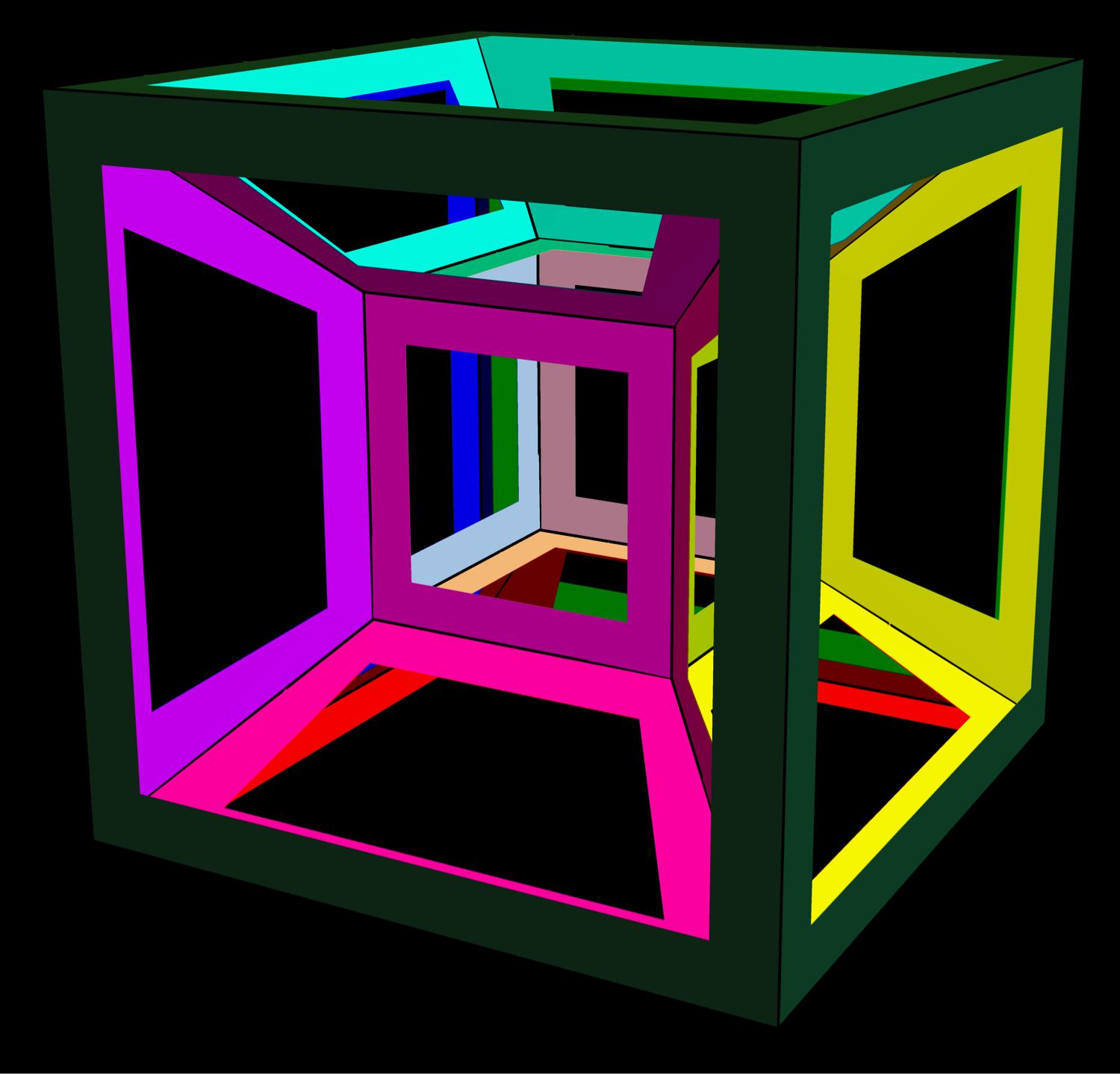}
\label{436_colouring_view_1}
}
\quad
\subfloat[Colouring of the hypercube, view 2.]
{
\includegraphics[width=0.3\textwidth]{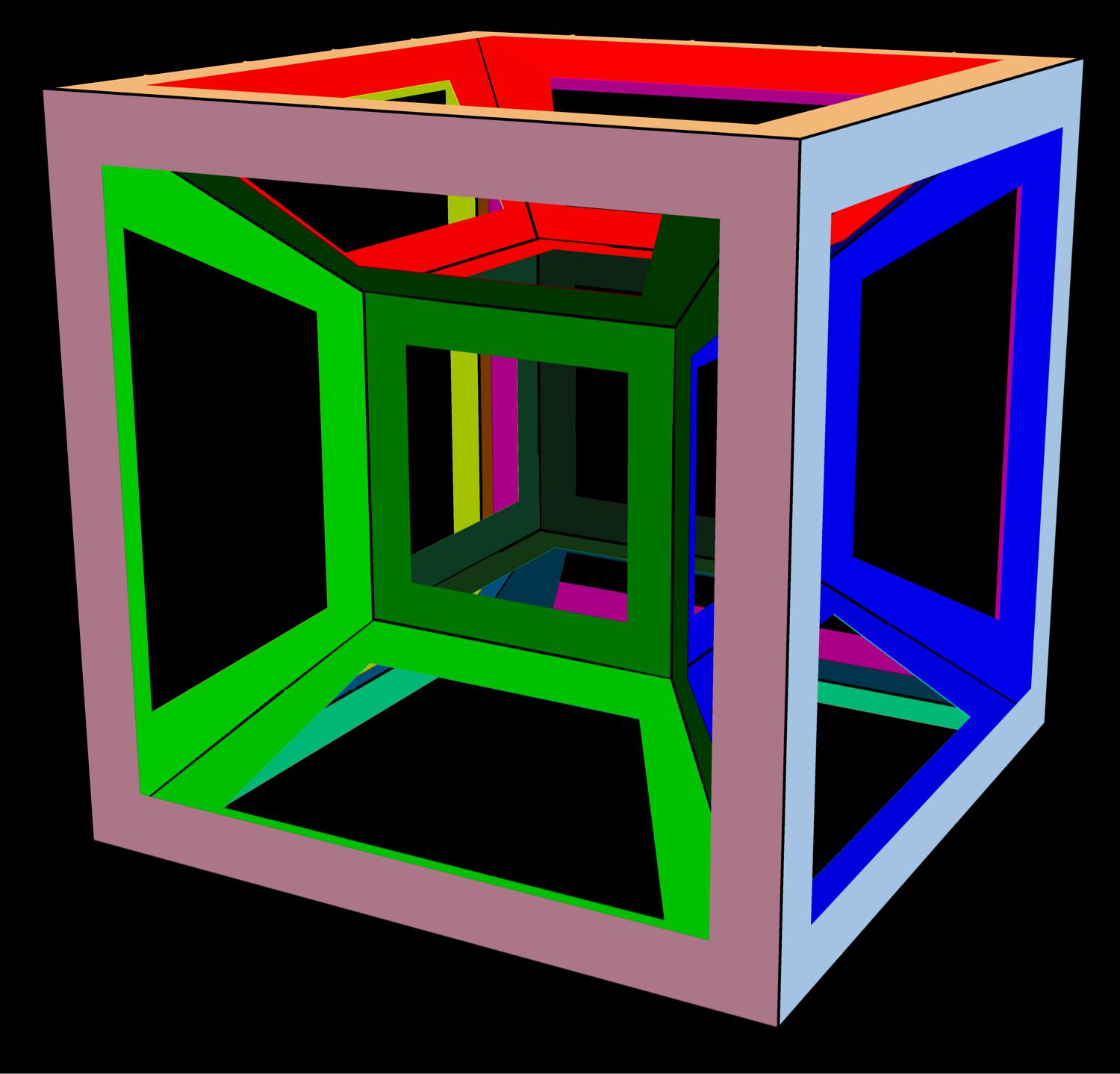}
\label{436_colouring_view_2}
}
\caption{Any colouring of the hypercube can be mapped onto a colouring of the $\{4,3,6\}$ hyperbolic honeycomb.}
\label{436_colouring}
\end{wrapfigure}

In Figure \ref{436}, we colour the cells using eight colours, in an interesting pattern very special to the $\{4,3,6\}$ honeycomb. This comes from the observation that $\{4,3,6\}$ is a kind of branched cover of the $\{4,3,3\}$ honeycomb, in which three cubes are arranged around each edge. The honeycomb $\{4,3,3\}$ does not tile hyperbolic space; rather it is a honeycomb that tiles spherical space: it is the same as the honeycomb we get by radially projecting the cubical cells of the hypercube onto a circumscribing three-sphere in four-dimensional 
space. To be more precise, there is a continuous map, $F$ say, from $\{4,3,6\}$ to $\{4,3,3\}$, that maps each cube of $\{4,3,6\}$ to one of the eight cubes of $\{4,3,3\}$. We assign a different colour to each of the eight cubes of $\{4,3,3\}$, as in Figure \ref{436_colouring}, then colour each cube $c$ of $\{4,3,6\}$ by the colour of $F(c)$. 
Patterns in the colouring can be seen in Figures \ref{436_trunc_cubes} and \ref{436_trunc_cubes_tris_only}: first that cubes opposite each other around an edge have the same colour, and second that going in a straight line, from face to opposite face of each cube, we get back to the same colour after four cubes.

\begin{figure}[htb]
\centering
\includegraphics[width=0.8\textwidth]{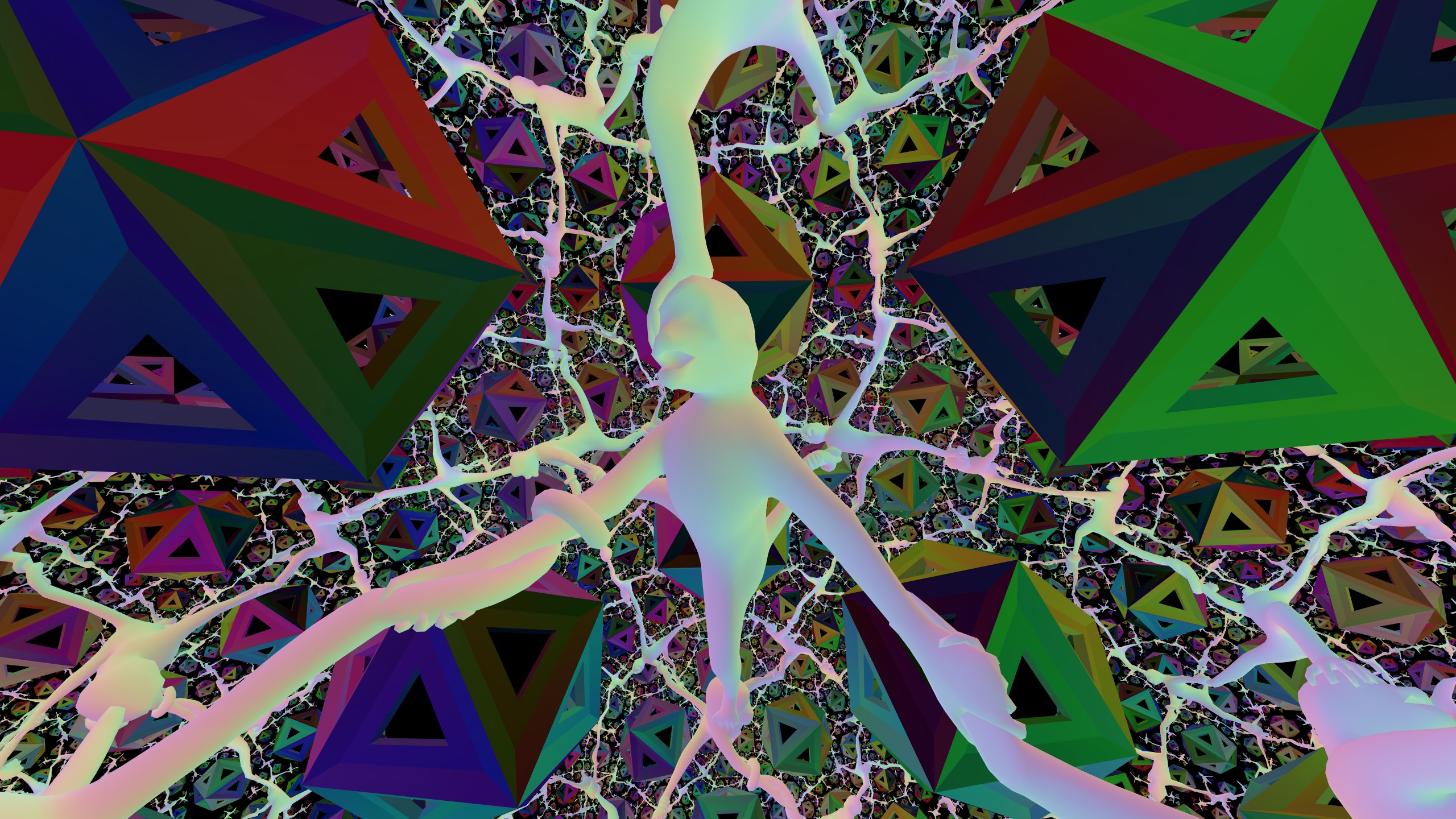}
\caption{A monkey in each cube of the $\{4,3,6\}$ honeycomb. Note the ring of six monkeys connected together around each edge of the honeycomb.}
\label{436_monkeys}
\end{figure}

Any pattern drawn on the hypercube can be lifted to the $\{4,3,6\}$ honeycomb. For example, our sculpture, \emph{More fun than a hypercube of monkeys}~\cite{bridges2014:143}, puts a monkey in each cubical cell of the hypercube. The lift of this sculpture is shown in Figure \ref{436_monkeys}.

\section{Virtual reality, parallel transport and the Levi-Civita connection}

\begin{wrapfigure}[12]{l}{0.3\textwidth}
  \vspace{-18pt}
  \centering
  \includegraphics[width=0.25\textwidth]{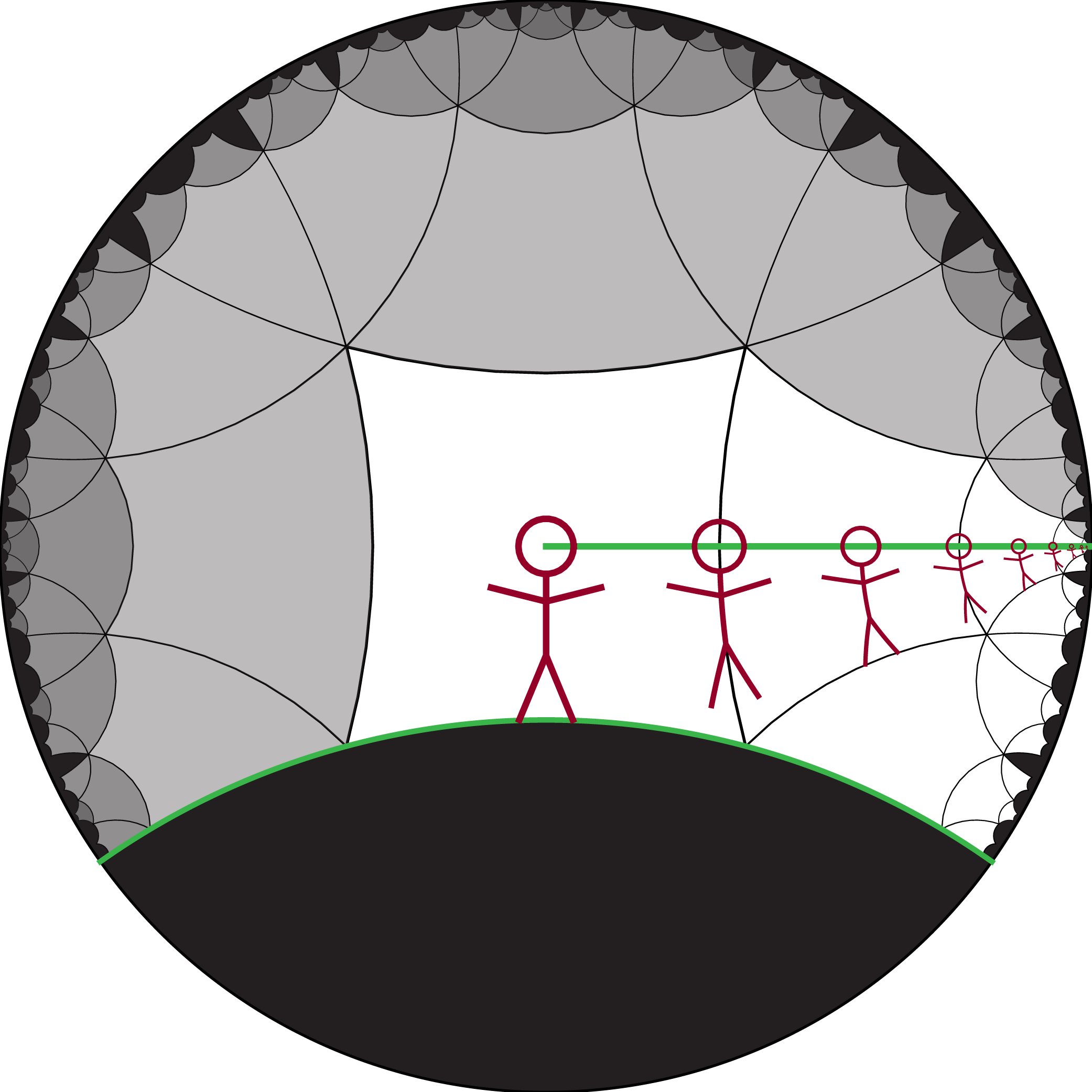}
  \caption{The floor falls out from under your feet as you travel along a geodesic.} 
  \label{Fig:floor_falling}
\end{wrapfigure}

The physicality of a virtual reality system with positional tracking gives us a visceral sense of some otherwise abstract phenomena. In a curved space, for example, two neighbouring geodesics that start with parallel velocities (tangent vectors) end up diverging if the space is negatively curved. Suppose that in the simulation, the user is standing on a floor consisting of a geodesic plane in $\HH^3$. When they walk forward in real-life, in the simulation their head follows a geodesic that starts out with velocity parallel to the floor, and which therefore diverges from the floor. This leads to the sensation that the floor is falling out from under your feet. See Figure \ref{Fig:floor_falling}.

This phenomenon is a consequence of \emph{parallel transport} -- as a vector is moved through space along a curve it stays parallel to itself and has constant magnitude. A formal definition of geodesics is that they are curves that parallel transport their own tangent vectors. How might we go about constructing geodesic from this notion? To move along a manifold in a path in a given direction, we must know how the velocity changes as we move parallel to the path. On a differentiable manifold $M$ with metric $g$, this is formalised by the notion of the Levi-Civita connection $\nabla_X$, which is the unique \emph{covariant derivative} -- the derivative in the manifold in the direction of $X$, a vector in the tangent space of $M$ -- that preserves the metric $\nabla g=0$ and is torsion free. 

\begin{wrapfigure}[12]{r}{0.45\textwidth}
  \vspace{-30pt}
  \centering
\subfloat[Parallel transport in $\EE^2$.]
{
\includegraphics[width=0.15\textwidth]{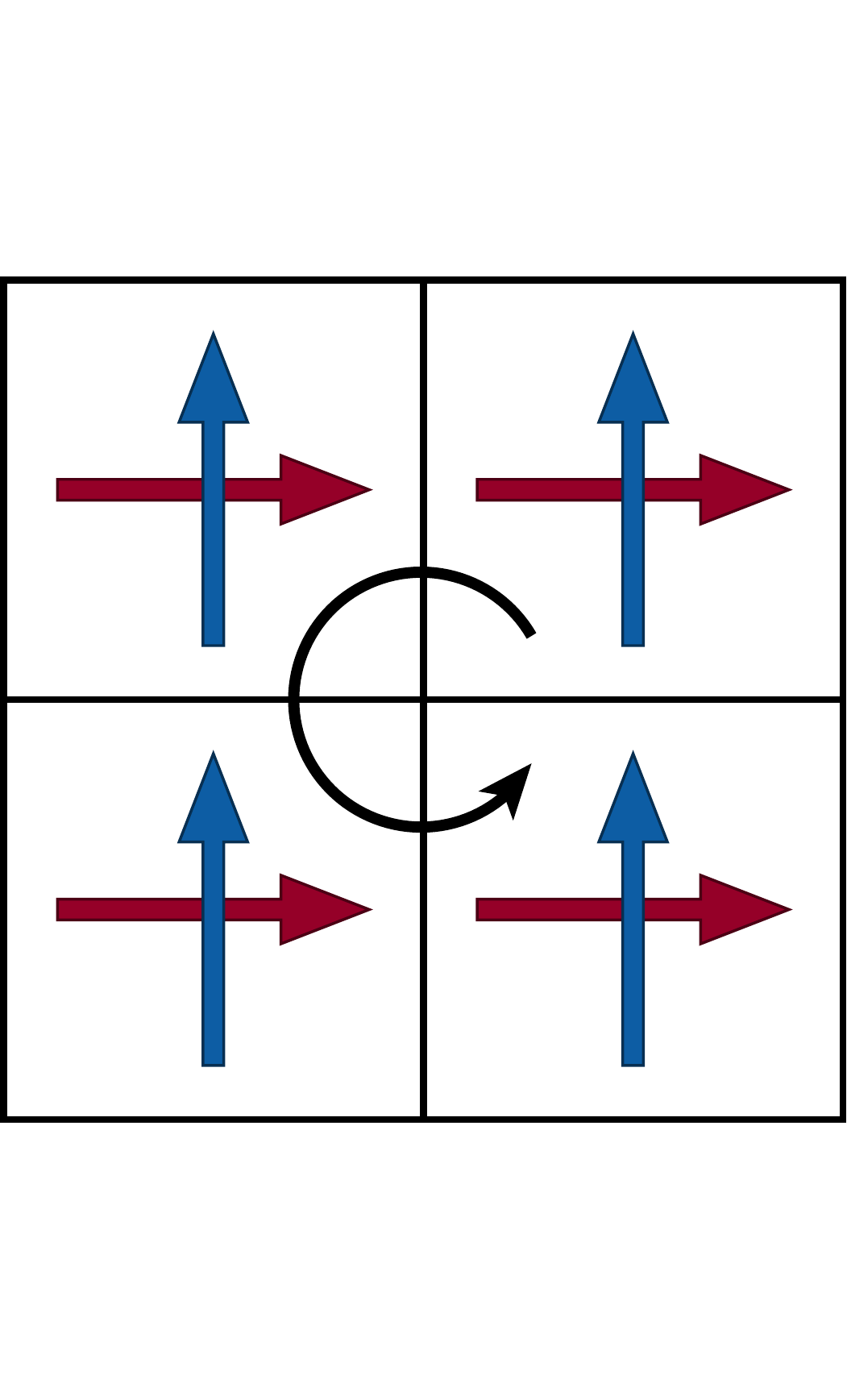}
\label{euclidean_connection}
}
\quad
 \subfloat[Parallel transport rotates the frame in $\HH^2$, shown in the Poincar\'e disk model.]
{
\includegraphics[width=0.25\textwidth]{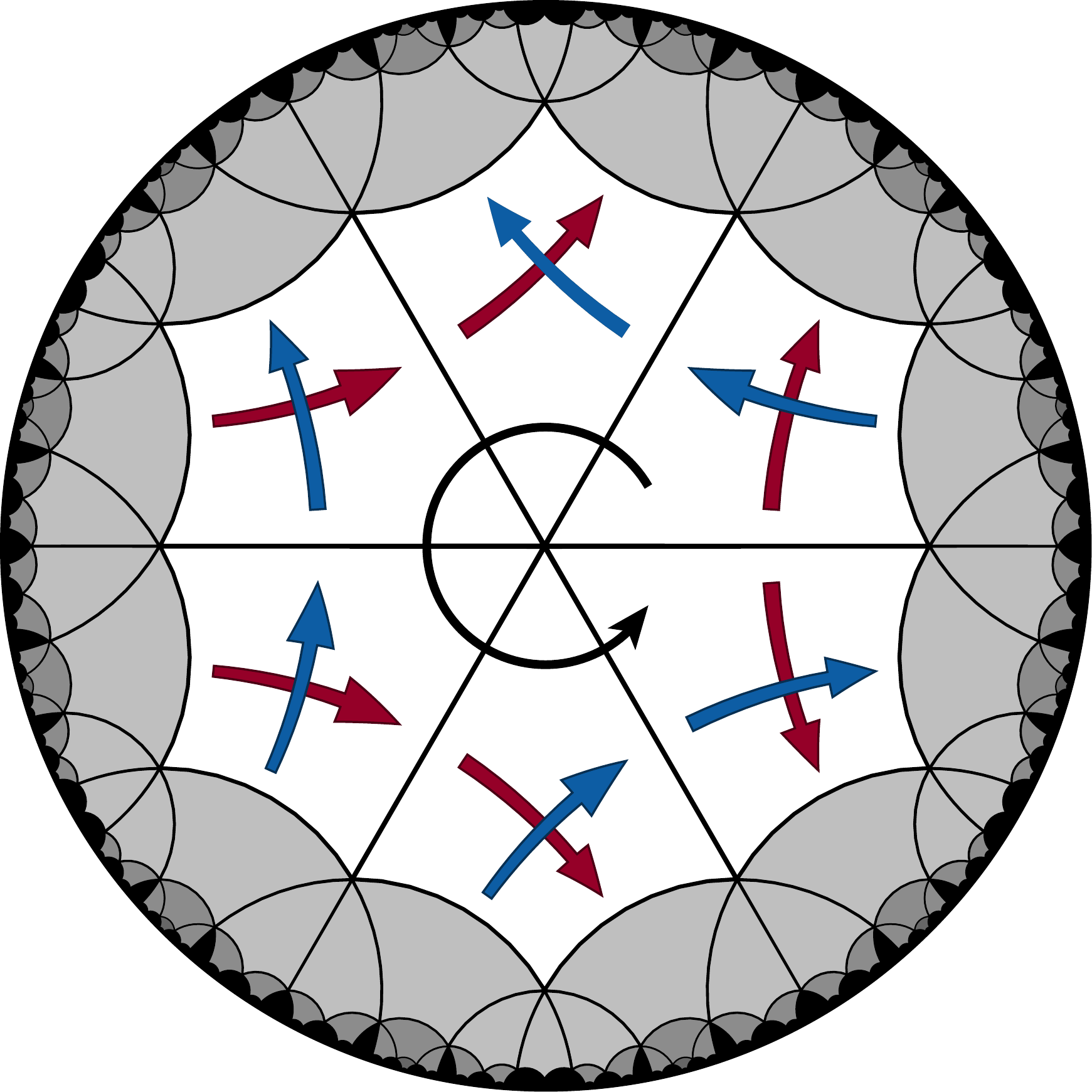}
\label{h2_connection}
}
  \caption{Walking around an edge.} 
  \label{Fig:metric_connection}
\end{wrapfigure}
\begin{wrapfigure}[45]{l}{0.4\textwidth}
\vspace{-10pt}
\centering
\subfloat[$\HH^3$ initial view.]
{
\includegraphics[width=0.35\textwidth]{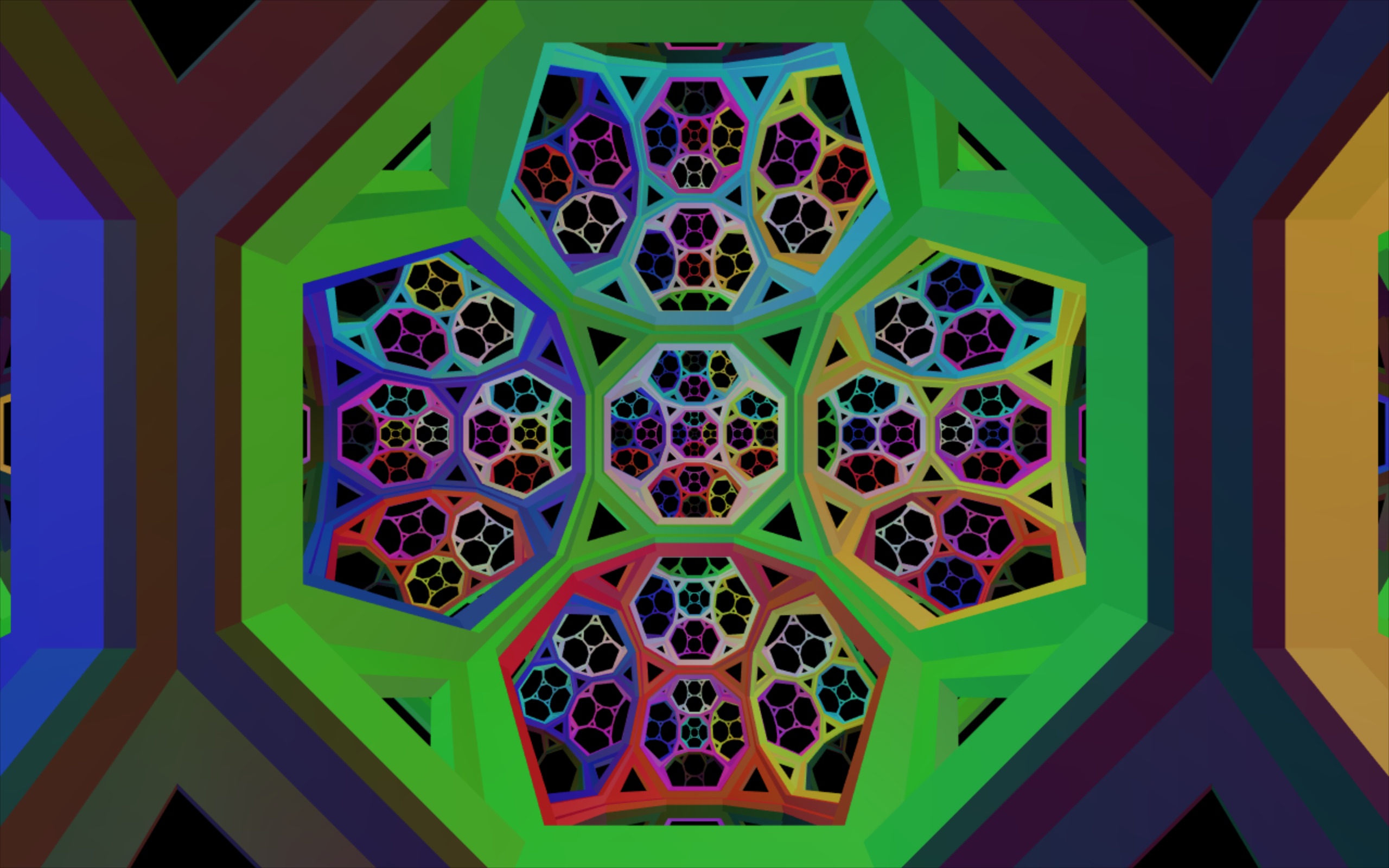}
\label{H3_parallel_transport_1}
}

\vspace{-7pt}
\subfloat[$\HH^3$ after moving right 0.5.]
{
\includegraphics[width=0.35\textwidth]{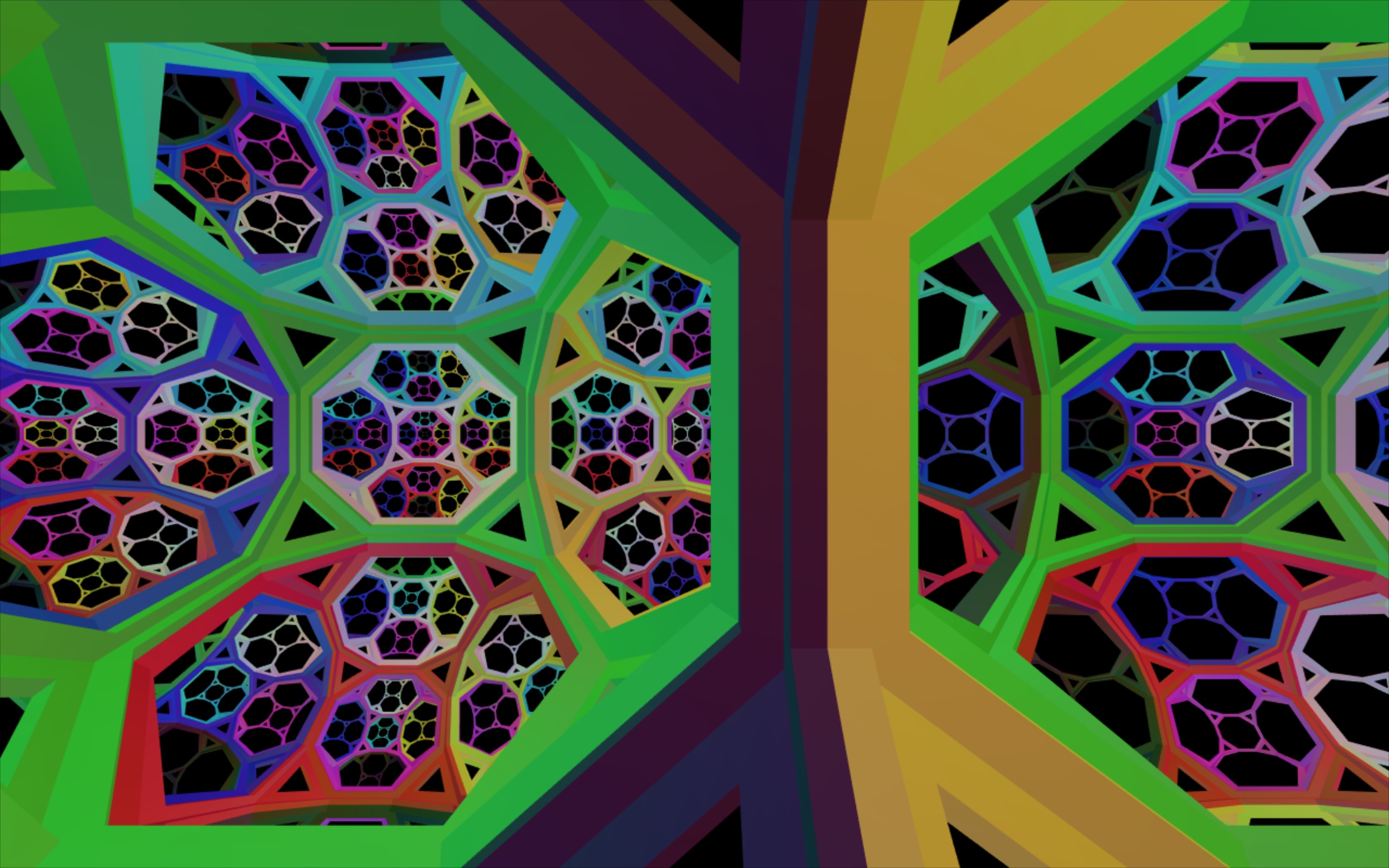}
\label{H3_parallel_transport_2}
}

\vspace{-7pt}
\subfloat[$\HH^3$ after moving up 0.5.]
{
\includegraphics[width=0.35\textwidth]{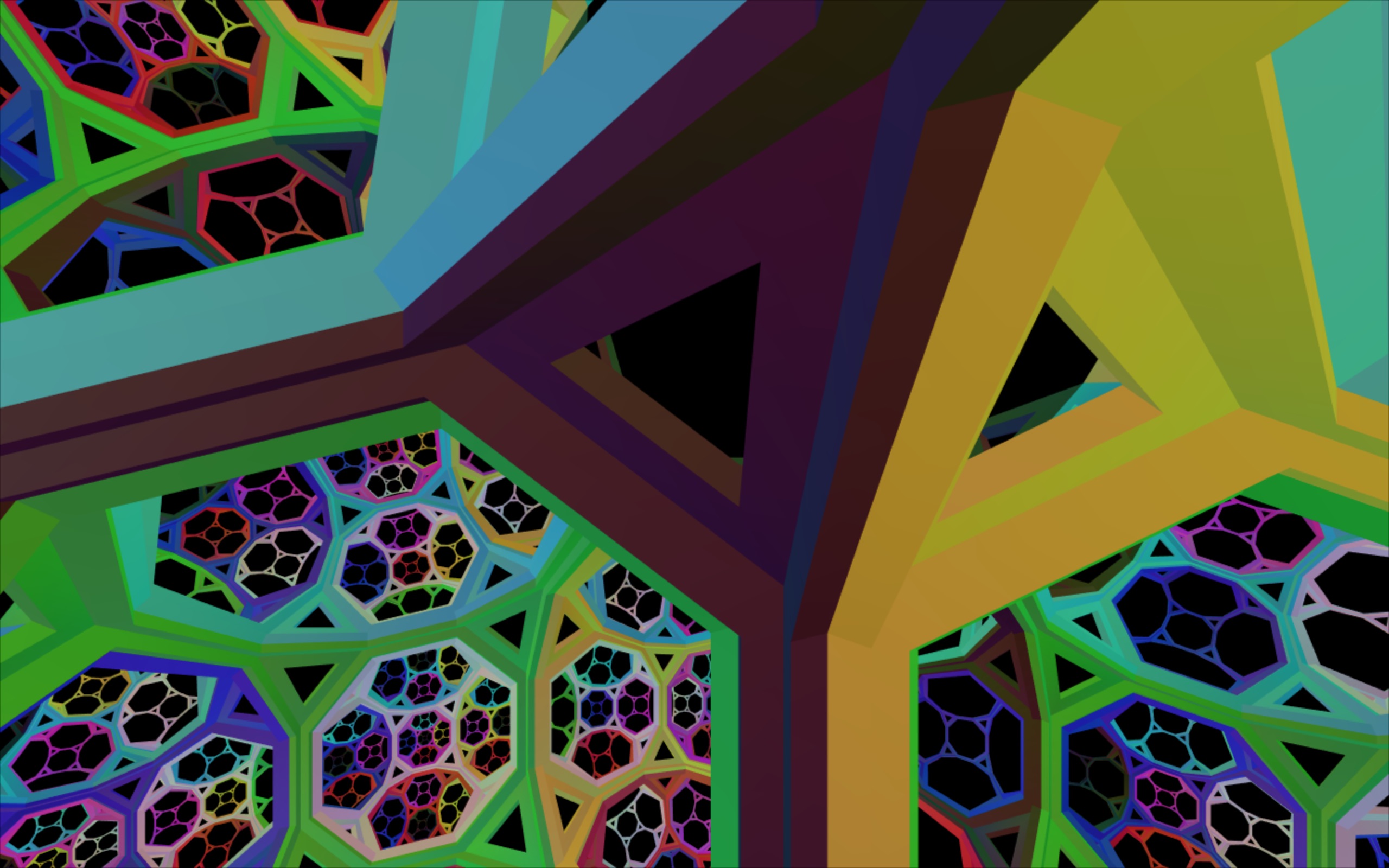}
\label{H3_parallel_transport_3}
}

\vspace{-7pt}
\subfloat[$\HH^3$ after moving left 0.5.]
{
\includegraphics[width=0.35\textwidth]{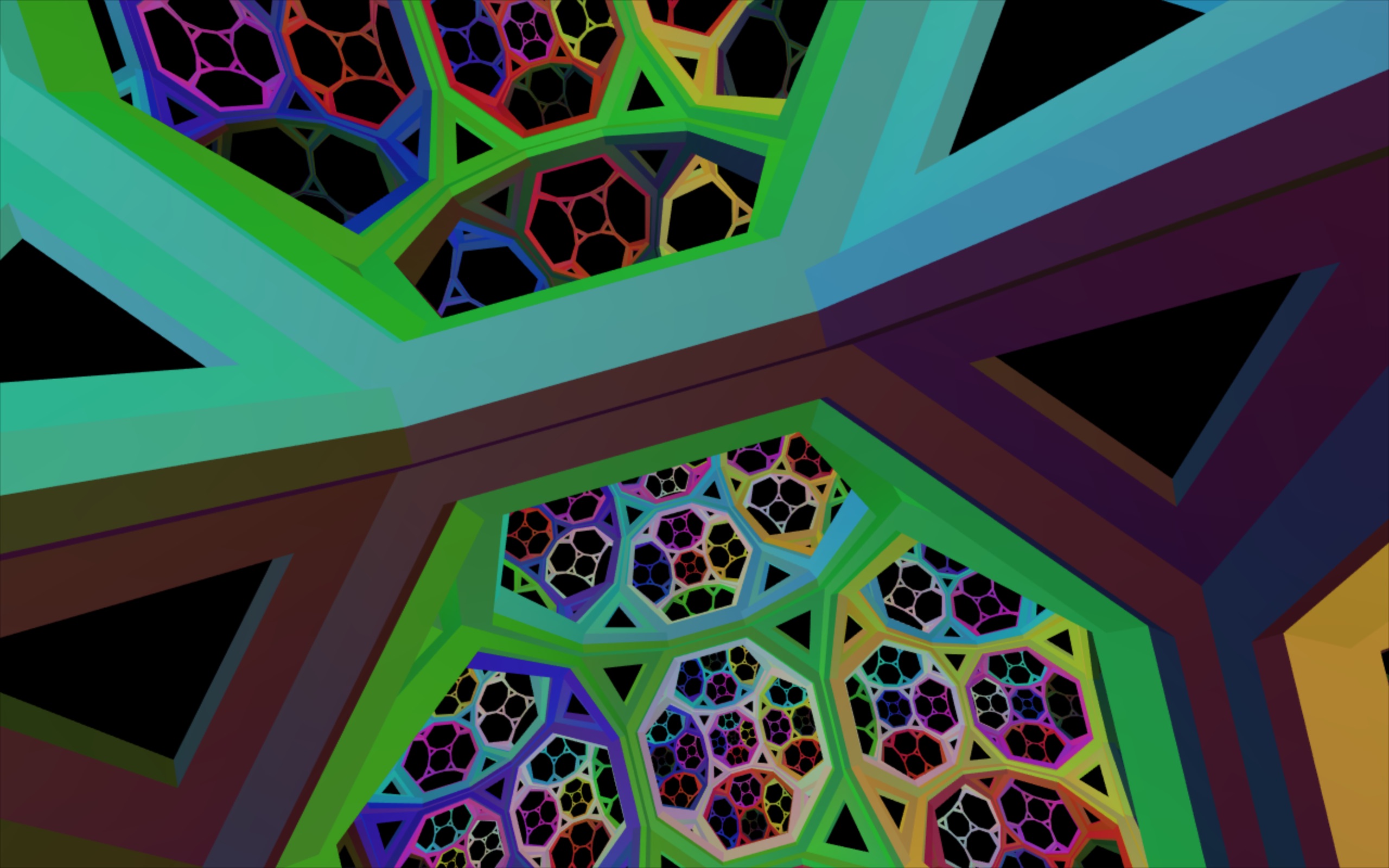}
\label{H3_parallel_transport_4}
}

\vspace{-7pt}
\subfloat[$\HH^3$ after moving down 0.5.]
{
\includegraphics[width=0.35\textwidth]{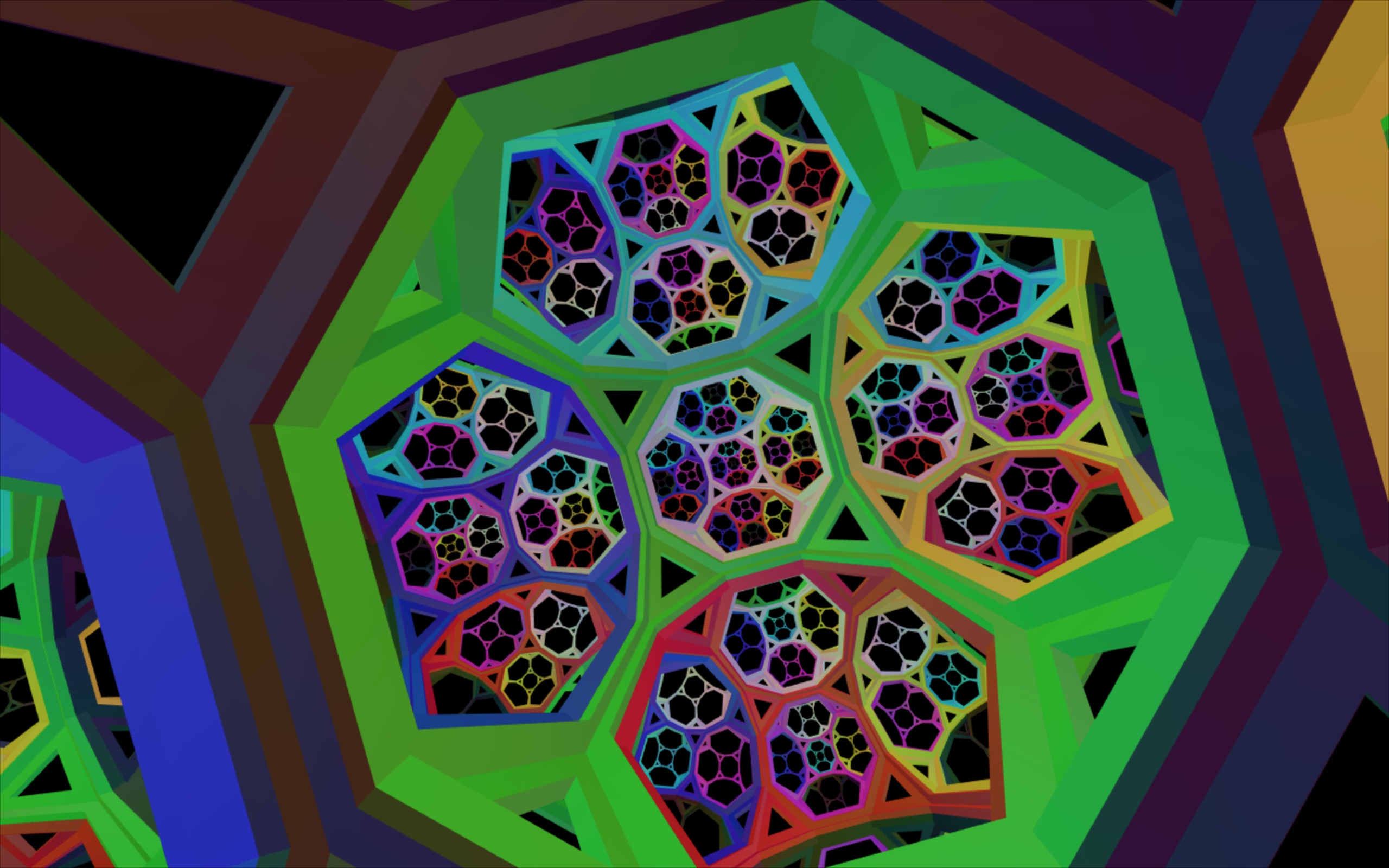}
\label{H3_parallel_transport_5}
}
\caption{Parallel transport rotates reference frames in curved space.}
\label{Fig:parallel_transport}
\end{wrapfigure}

Another, unexpected phenomenon we encounter in the virtual reality experience stems from parallel transport of other vectors. When we experience the world, we are aware not only of the vector that points in the direction we are looking but also of the vectors that point up and down, left and right. Moving along a path in the virtual reality $\HH^3$ space, these vectors get transported as well. We have a fixed sense of which direction ``up" is, but this direction can rotate with respect to the world. See Figures \ref{Fig:metric_connection} and \ref{Fig:parallel_transport}. In particular, this means that certain movements in $\EE^3$ produce a rotation of 
the floor of a room drawn in $\HH^3$, so that it no longer appears to coincide with the real-life floor the user is walking on.

These phenomena make $\HH^3$ a somewhat confusing place to live in, at least as a visitor from $\EE^3$.
There may be ways to ``hack'' the simulation to solve the problems of the virtual floor falling away or rotating away
from the real-life floor. To ``fix'' the angle of the floor changing, we could artificially rotate the virtual view so that the orientation of the virtual camera relative to the virtual floor always agrees with the orientation of the headset relative to the real-life floor. 
Alternatively, we could avoid both problems by tracking the point directly between the user's feet rather than their head as it moves through space, and for every frame offset the position of the
camera up from the feet to the head. 
These are both somewhat artificial fixes however, and would preclude the user from experiencing the effects of parallel transport. 

\newpage

$\phantom{aaa}$
\vspace{-43pt}

\section{Future directions}
In addition to putting more recognisable objects and architecture into our simulations and allowing users to interact with objects, we would like to build similar simulations for the other Thurston geometries. Each of these geometries presents a unique challenge. Both $S^3$ and $S^2\times \EE$ have a multivalued exponential maps, thus we need to calculate the exponential map both in front of the viewer and behind them in order to draw a complete image on the screen. Nil, Solv and $\widetilde{\text{PSL}_2(\RR)}$ don't have the ubiquitous standard models that spherical and hyperbolic space have. In our future implementations of Thurston geometries, we will use models of Emil Moln{\'a}r \cite{molnar_geometries} to create the exponential map as well as the set of isometries.

The most natural extension of the work on $\HH^3$ is to the product space $\HH^2\times \EE$ -- the cartesian product of the hyperbolic plane with the euclidean line. 
We discuss our simulation of $\HH^2\times\EE$ in our second paper in this volume~\cite{our_paper_h2xe}.

\vspace{-2pt}

\bibliographystyle{plain}
\bibliography{h2xr.bib}

\begin{thebibliography}{1}

\bibitem{our_paper_h2xe}
Vi~Hart, Andrea Hawksley, Elisabetta~A. Matsumoto, and Henry Segerman.
\newblock Non-euclidean virtual reality \textrm{II}: explorations of
  $\mathbb{H}^2\times\mathbb{E}$.
\newblock In {\em Proceedings of Bridges 2017: Mathematics, Music, Art,
  Architecture, Culture}. Tessellations Publishing, 2017.

\bibitem{bridges2014:143}
Vi~Hart and Henry Segerman.
\newblock The quaternion group as a symmetry group.
\newblock In Gary Greenfield, George Hart, and Reza Sarhangi, editors, {\em
  Proceedings of Bridges 2014: Mathematics, Music, Art, Architecture, Culture},
  pages 143--150. Tessellations Publishing, 2014.
\newblock \url{http://archive.bridgesmathart.org/2014/bridges2014-143.html}.

\bibitem{molnar_geometries}
Emil Moln{\'a}r.
\newblock The projective interpretation of the eight 3-dimensional homogeneous
  geometries.
\newblock {\em Beitrage zur Algebra und Geometrie (Contributions to Algebra and
  Geometry)}, 38(2):261--288, 1997.

\bibitem{visualizing_hyperbolic_honeycombs}
Roice Nelson and Henry Segerman.
\newblock Visualizing hyperbolic honeycombs.
\newblock arXiv:1511.02851.

\bibitem{perelman1}
Grisha Perelman.
\newblock The entropy formula for the {R}icci flow and its geometric
  applications.
\newblock {arXiv:0211159}.

\bibitem{thurston_book}
William~P. Thurston.
\newblock {\em Three-Dimensional Geometry and Topology}.
\newblock Princeton Univ. Press, 1997.

\bibitem{curved_spaces}
Jeff Weeks.
\newblock {Curved Spaces}.
\newblock a flight simulator for multiconnected universes, available from
  \url{http://www.geometrygames.org/CurvedSpaces/}.

\bibitem{weeks_real-time_rendering}
Jeff Weeks.
\newblock Real-time rendering in curved spaces.
\newblock {\em IEEE Computer Graphics and Applications}, 22(6):90--99, 2002.

\end{thebibliography}

\end{document}